\documentclass[10pt]{article}

\usepackage{epsfig,amsmath}

\usepackage{amssymb}

\newcommand{\proof}{\noindent{\sl Proof.\/}\ }
\newcommand{\finprf}{\unskip\null\hfill$\square$\vskip 0.3cm}
\newcommand{\un}{{\rm 1\kern -2.5pt l}}

\newtheorem{theorem}{\textsc{Theorem}}[section]

\newtheorem{corollary}[theorem]{\textsc{Corollary}}

\begin{document}
\title{{\bf \Large On the one-dimensional parabolic obstacle problem\\
    with variable coefficients}}
\author{\normalsize\textsc{A. Blanchet}\footnotemark[1]{ }\footnotemark[2],
\normalsize\textsc{J. Dolbeault}\footnotemark[1],
\normalsize\textsc{R. Monneau}\footnotemark[2]}
\date{}
\maketitle
\footnotetext[1]{\it CEREMADE, Universit\'e Paris Dauphine, place de Lattre de Tassigny, 75775 Paris C\'edex 16, France}
\footnotetext[2]{\it CERMICS, Ecole Nationale des Ponts et Chauss\'ees, 6 et 8 avenue
    Blaise Pascal,
  Cit\'e Descartes Champs-sur-Marne, 77455
 Marne-la-Vall\'ee  C\'edex 2, France}
\centerline{\small{\bf{Abstract}}}
\noindent{\small{This note is devoted to continuity results of the time
    derivative of the solution to the one-dimensional parabolic obstacle problem with
    variable coefficients. It applies to the smooth fit principle in
    numerical analysis and in financial mathematics. It relies on
    various tools for the study of free boundary problems: blow-up method, monotonicity
    formulae, Liouville's results.}}\hfill\break

\noindent{\small{\bf{AMS Classification:}}} {\small{35R35.}}\hfill\break
\noindent{\small{\bf{Keywords:}}} {\small{parabolic obstacle problem,
    free boundary, blow-up, Liouville's result, monotoni\-city formula,
    smooth fit.}}\hfill\break

\section{Introduction}
Consider a parabolic obstacle problem in an open set. We look for local properties, which do not depend on the
boundary conditions and the initial conditions, but only depend on the
equation in the interior of the domain. Consider a function $u$
with a one-dimensional space variable~$x$ in $Q_1(0)$ where by $Q_r(P_0)$ we denote the
parabolic box of radius $r$ and of centre $P_0=(x_0,t_0)$:
$$Q_r(P_0)=\left\{(x,t)\in \mathbb R^2,\quad |x-x_0|< r,\quad |t-t_0|< r^2\right\}\;.$$

Assume that $u$ is a solution of the one-dimensional parabolic obstacle
problem with variable coefficients:
\begin{equation}\label{eq:1}
\left\{
\begin{array}{l}
a(x,t)u_{xx}+b(x,t)u_x +c(x,t)u -u_t = f(x,t) \cdot
\un_{\left\{u>0\right\}} \quad \mbox{a.e. in } Q_1(0)\vspace{.1cm}\\
u \ge 0  \quad \mbox{a.e. in } Q_1(0)
\end{array}
\right.
\end{equation}
where $u_t,u_x,u_{xx}$ respectively stand for $\frac{\partial u}{\partial
  t}$, $\frac{\partial u}{\partial
  x}$, $\frac{\partial^2 u}{\partial
  x^2}$, and $\un_{\left\{u>0\right\}}$ is the characteristic function of the positive set of $u$.
Here the free boundary $\Gamma$ is defined by
$$\Gamma = \left(\partial \left\{u=0\right\}\right) \cap Q_1(0)\;.$$
To simplify the presentation, we assume that the coefficients
\begin{equation}\label{eq:2}
a,\;b,\;c \mbox{ and }f\mbox{ are }  C^1 \mbox{ in }(x,t) \;,  
\end{equation} 
but H\"{o}lder continuous would be sufficient in what follows.

A natural assumption is that the differential operator is uniformly
elliptic, {\it i.e.} the coefficient $a$ is bounded from below by
zero. If we do not make further assumptions on $a$ and on $f$, we cannot expect any good
property of the free boundary $\Gamma$. Suppose that:
\begin{equation}\label{eq:3}
\exists \,\delta >0,\quad a(x,t)\ge \delta,\; f(x,t) \ge \delta \;
\mbox{ a.e. in }Q_1(0) \;.  
\end{equation}
Up to a reduction of the size of the box (see \cite{BDM}), any weak solution $u$ of
(\ref{eq:1}) has a bounded first derivative in time and bounded first
and second
derivatives in space. Assume therefore that this property holds on the
initial box:
\begin{equation}\label{eq:4}
|u(x,t)|, \ \left|u_t(x,t)
\right|,\  
\left|u_x(x,t) \right| \mbox{ and }
\left|u_{xx}(x,t) \right| \mbox{ are bounded in } Q_1(0) \;.  
\end{equation}

This problem is a generalisation to the case of an operator with variable
coefficients of Stefan's problem (case where the parabolic operator is ${\partial^2}/{\partial
  x^2}-{\partial}/{\partial t}$). Stefan's problem  describes
the interface of ice and water (see \cite{MR81g:49013,rodrigues,friedlebouquin}). The problem with variable coefficients arises in the pricing of
american options in financial mathematics (see \cite{b&s,bensou,van,ll,jll,villeneuve,achdou}).

If $P$ is a point
such that $u(P)>0$, by standard parabolic
estimates $u_t$ is continuous in a neighbourhood
of $P$. On the other hand if $P$ is in the interior of the region
$\left\{u=0\right\}$, $u_t$ is obviously continuous. The only difficulty is therefore the regularity on the
free boundary~$\Gamma$. By assumption $u_t$ is bounded but may be {\it discontinuous} on $\Gamma$. The regularity of $u_t$ is a crucial question to apply the ``smooth-fit principle''
which amounts to the $C^1$ continuity of the solution at the
free boundary. This principle is often assumed, especially in the papers dealing with numerical analysis (see P. Dupuis and
H. Wang \cite{dupuiswang} for example).

In a recent work L. Caffarelli,
A. Petrosyan and H. Shagholian \cite{CPS} prove the
$\mathcal C^\infty$ regularity of the free boundary locally around some
points which are energetically characterised, without any sign
assumption neither on $u$ nor on its time derivative. This result holds in
higher dimension but in the case of constant coefficients. We use tools
similar to the ones of \cite{CPS} and the ones the last author developed previously for the
elliptic obstacle problem in \cite{mono2}. Our main result is the following:
\begin{theorem}[Continuity of $u_t$
    for almost every time]\label{th:1}
Under assumptions (\ref{eq:1})-(\ref{eq:2})-(\ref{eq:3})-(\ref{eq:4}), for
almost every time $t$, the function $u_t$ is
continuous on $Q_1(0)$. 
\end{theorem}
This result is new, even in the case of constant coefficients. The continuity of $u_t$ cannot be obtained everywhere in $t$, as shown by the following example. Let $u(x,t)=\max\{0,-t\}$. It satisfies $u_{xx}-u_t
=\un_{\left\{u>0\right\}}$ and its time derivative is obviously discontinuous at $t=0$.

If additionally we assume that $u_t\ge 0$ we achieve a more precise result:
\begin{theorem}[Continuity of $u_t$ for all $t$ when $u_t\ge 0$]\label{th:2}
Under assumptions (\ref{eq:1})-(\ref{eq:2})-(\ref{eq:3})-(\ref{eq:4}), if $u_t \ge 0$ in $Q_1(0)$  
then $u_t$ is continuous everywhere in $Q_1(0)$.
\end{theorem}
The assumption that $u_t\ge 0$ can be established in some special cases (special initial conditions, boundary conditions, and time independent coefficients). See for example the results of Friedman \cite{frie1d2}, for further results on the one-dimensional parabolic obstacle problem with particular initial conditions.

In Section \ref{blow-up} we introduce blow-up sequences, which are a kind
of zoom at a point of the free boundary. They converge, up to a
sub-sequence, to a solution on the whole space of the obstacle problem
with constant coefficients. Thanks to a monotonicity formula for an
energy we prove in Section \ref{fmlg} that the blow-up limit
is scale-invariant. This allows us to
classify in Section \ref{liouville} all possible blow-up limits in a Liouville's
theorem. Then we sketch the proof of Theorem \ref{th:1}. We
even classify energetically the points of the free boundary into the set
of regular and singular points. In Section \ref{singulier} we prove the
uniqueness of the blow-up limit at singular points. Then we give the
sketch of the proof of the Theorem \ref{th:2}. For further details we refer to \cite{BDM}.
\section{The notion of blow-up}\label{blow-up}
\begin{figure}[h!]
   \begin{center}
      \includegraphics[width=8cm]{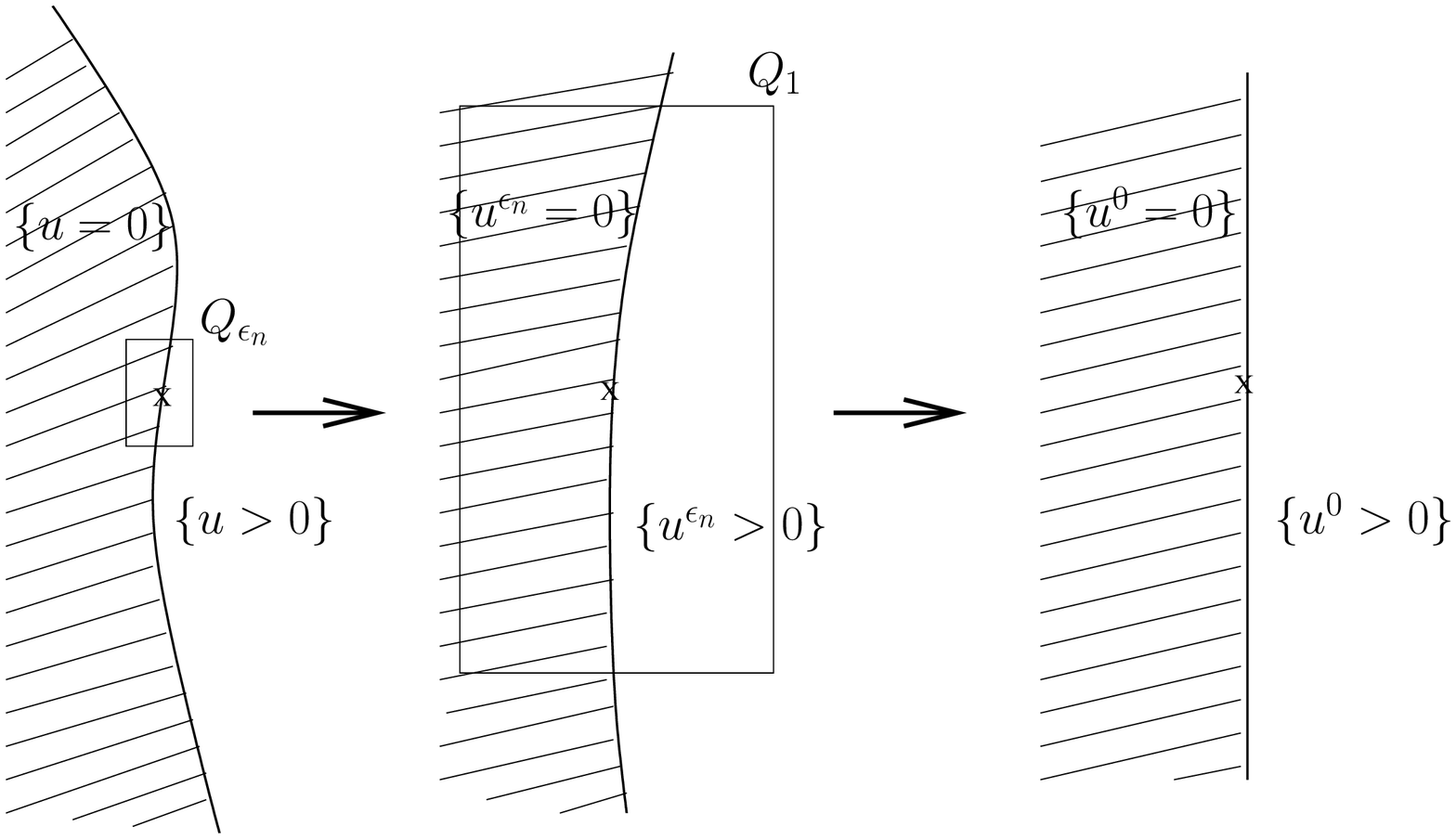}
   \end{center}
   \caption{\footnotesize Blow-up}\label{fig.1}
\end{figure}

Given a point $P_0=(x_0,t_0)$ on the free boundary $\Gamma$, we can define
the blow-up sequence by
\begin{equation}\label{eq:5}
u^\varepsilon_{P_0}(x,t)=\frac{u(x_0 +\varepsilon x, t_0 +\varepsilon^2
  t)}{\varepsilon^2}\;, \varepsilon>0 \;.  
\end{equation}
Roughly speaking the action of this rescaling is to zoom on the free boundary at
scale~$\varepsilon$ (see figure \ref{fig.1}).

By assumption, $u(P_0)=0$. Because $u$ is non-negative, we also have
  $u_x(P_0)=0$. Moreover $u^\varepsilon_{P_0}$ has a bounded first derivative
  in time and bounded second derivatives in space. For this reason, using
  Ascoli-Arzel\`a's theorem, we can find a sequence $(\varepsilon_n)_n$ which converges to zero such that $\left(u^{\varepsilon_n}_{P_0}\right)_{n}$ converges on every compact set of
  $\mathbb R^2= \mathbb R_x\times \mathbb R_t$ to a function $u^0$ (called the blow-up limit)
  and which {\it a priori} depends on the choice of the sequence
  $(\varepsilon_n)_{n}$.

The limit function $u^0$ satisfies the parabolic obstacle problem with constant coefficients on the whole
  space-time:
$$a(P_0)u^0_{xx}-u^0_t = f(P_0)\cdot \un_{\left\{u>0\right\}} \quad
\mbox{in}\quad \mathbb R^2\;.$$
By the non-degeneracy assumption (\ref{eq:3}), it is possible to
prove that $0\in \partial\left\{u^0=0\right\}$.

To characterise the blow-up limit $u^0$, we need to come back to the original equation satisfied by $u$ and to
obtain additional estimates. In order to simplify the presentation we make a
much stronger assumption on $u$: assume that $u$ is a solution on the whole space-time of the equation with constant
coefficients $a\equiv1$, $ f\equiv 1$, $b\equiv 0$ and $c\equiv0$:
\begin{equation}\label{etoilejean}
u_{xx}-u_t = \un_{\left\{u>0\right\}} \quad
\mbox{in}\quad \mathbb R^2\;.
\end{equation}
Without this assumption, all tools have to be localised. See \cite{BDM} for more details.
\section{A monotonicity formula for energy}\label{fmlg}
For every time $t<0$, we define the quantity
$${\mathcal E}(t;u)=\int_{\mathbb R}\left\{\frac{1}{-t}\left(|u_x(x,t)|^2 + 2\,
    u(x,t)\right) -\frac{1}{t^2}\,u^2(x,t)\right\}G(x,t)\ dx$$
where $G$ satisfies the backward heat equation $G_{xx}+G_t = 0$ in
    $\left\{t<0\right\}$ and is given by
$$G(x,t)=\frac{1}{2\sqrt{\pi (-t)}}\exp \left(\frac{- x^2}{4(-t)}\right)\;.$$
\begin{theorem}[Monotonicity formula for energy]\label{th:3}
Assume that $u$ is a solution of~(\ref{etoilejean}). The function ${\mathcal E}$ is non-increasing in time for $t<0$, and
satisfies
\begin{equation}\label{eq:10}
\frac{d}{dt}{\mathcal E}(t;u) = -\frac{1}{2(-t)^3}\int_{\mathbb R} |{\mathcal
  L}u(x,t)|^2 G(x,t) \ dx   
\end{equation}
where
$${\mathcal
  L}u(x,t) = -2\,u(x,t) + x\cdot u_x(x,t) +2\,t\cdot  u_t(x,t)\;.$$
\end{theorem}
A similar but different energy is introduced in \cite{CPS,WEISS99}.
\begin{corollary}[Homogeneity of the blow-up limit]\label{cor:4}
Any blow-up limit $u^0$ of $(u^{\varepsilon_n}_{P_0})_{n}$ defined in (\ref{eq:5}), 
satisfies 
\begin{equation}\label{eq:6}
u^0(\lambda x,\lambda^2 t)= \lambda^2
u^0(x,t)\quad \mbox{for every}\quad x\in \mathbb R, \ t<0, \ \lambda >0 \;. 
\end{equation}
\end{corollary}
\proof
We prove it in the case $P_0=0$. The crucial property is the scale-invariance of $\mathcal E$:
$${\mathcal E}(\varepsilon_n^2 t;u)= {\mathcal E}(t;u^{\varepsilon_n}_0)\;.$$
Taking the limit $\varepsilon_n \to 0$, we get
\begin{equation}\label{eq:8}
{\mathcal E}(0^-;u):=\lim_{\tau \to 0\; \tau<0}{\mathcal E}(\tau;u)= {\mathcal E}(t;u^0) \quad \mbox{for every}\quad
t<0\;.
\end{equation}
From the monotonicity formula (\ref{eq:10}), we get ${\mathcal
  L}u^0(x,t) =0$ for $t <0$. This implies the homogeneity (\ref{eq:6})
  of $u^0$ for $t<0$.
\finprf
\section{A Liouville's theorem and consequences}\label{liouville}
Let
$$v_+(x,t) =\frac{1}{2} (\max\{0,x\})^2\;,$$
$$v_-(x,t) =\frac{1}{2} (\max\{0,-x\})^2\;,$$
and for $m\in [-1,0]$
$$v_m(x,t) = \left\{\begin{array}{ll}
\displaystyle{m\,t+\frac{1+m}{2}\,x^2} & \quad \mbox{if}\quad t\leq 0\;,\vspace{.3cm}\\
\displaystyle t\, V_m\left(\frac{|x|}{t}\right)>0 & \quad \mbox{if}\quad 0< t< C_m \cdot x^2\;,\vspace{.3cm}\\
0 & \quad \mbox{if}\quad t\geq  C_m \cdot x^2\;,
\end{array}\right.$$
where the coefficient $C_m$ is an increasing function of $m$, satisfying
$C_m =0$ if $m=-1$, and $C_m =+\infty $, if $m=0$. The precise expression
of $V_m$ is given in \cite{BDM}. In particular we get
$v_{-1}(x,t)=\max\{0,-t\}$ and $v_{0}(x,t)=\frac12 x^2$.
\begin{theorem}[Classification of global homogeneous solutions in $\mathbb R^2$]\label{classirr}
Let \hbox{$u^0\not\equiv 0$} be a non-negative solution of (\ref{etoilejean}) satisfying the homogeneity condition (\ref{eq:6}). Then $u^0$ is one of $v_+$, $v_-$ or $v_m$ for some $m\in [-1,0]$.
\end{theorem}

\begin{figure}[h!]
\begin{center}
\begin{minipage}[b]{0.17\linewidth}
\includegraphics[scale=0.3]{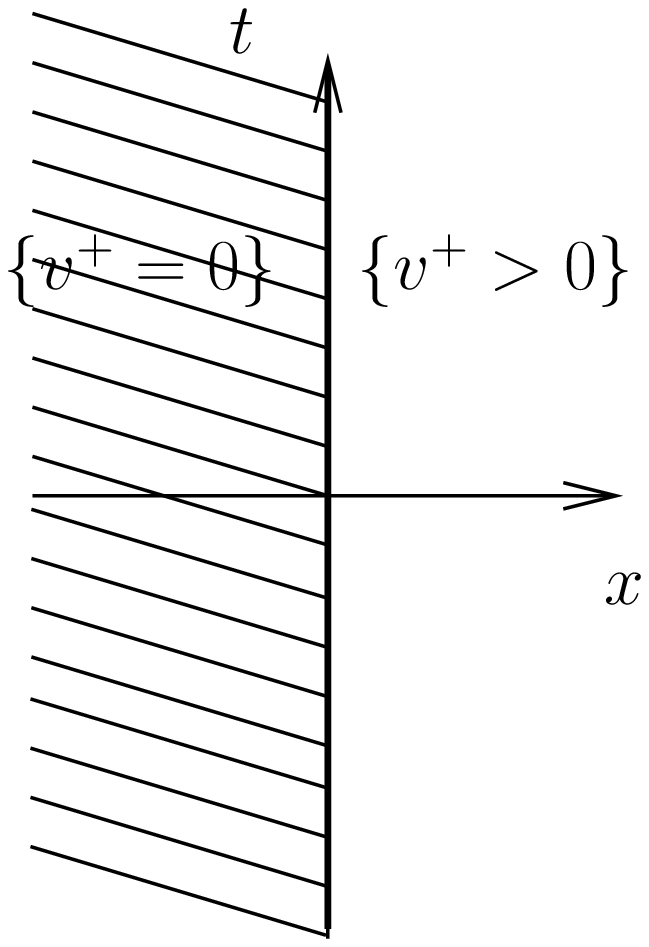}
\end{minipage}
\begin{minipage}[b]{0.17\linewidth}
\includegraphics[scale=0.3]{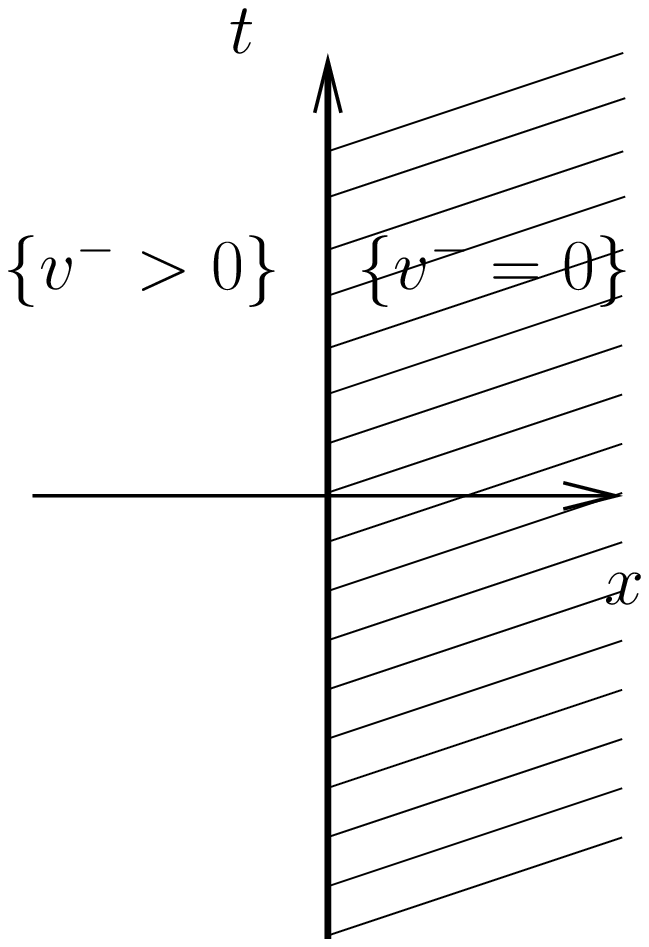}
\end{minipage}
\begin{minipage}[b]{0.22\linewidth}
\includegraphics[scale=0.3]{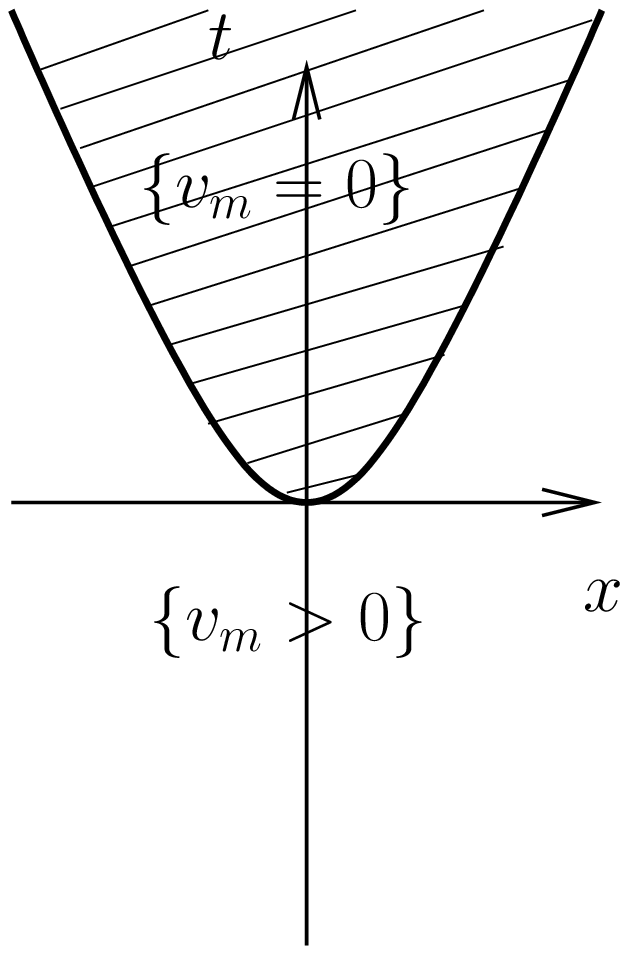}
\end{minipage}
\begin{minipage}[b]{0.17\linewidth}
\includegraphics[scale=0.3]{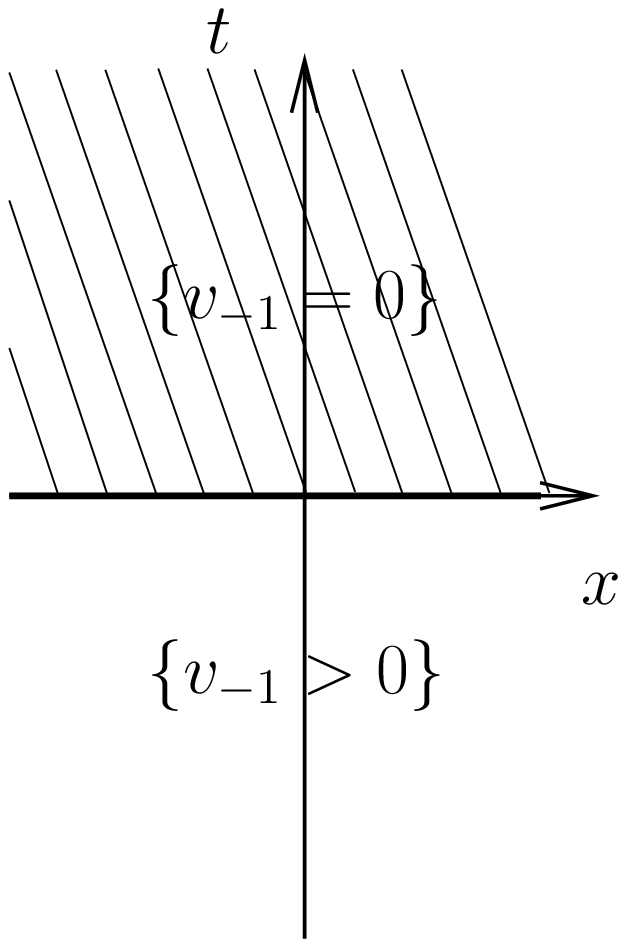}
\end{minipage}
\begin{minipage}[b]{0.17\linewidth}
\includegraphics[scale=0.3]{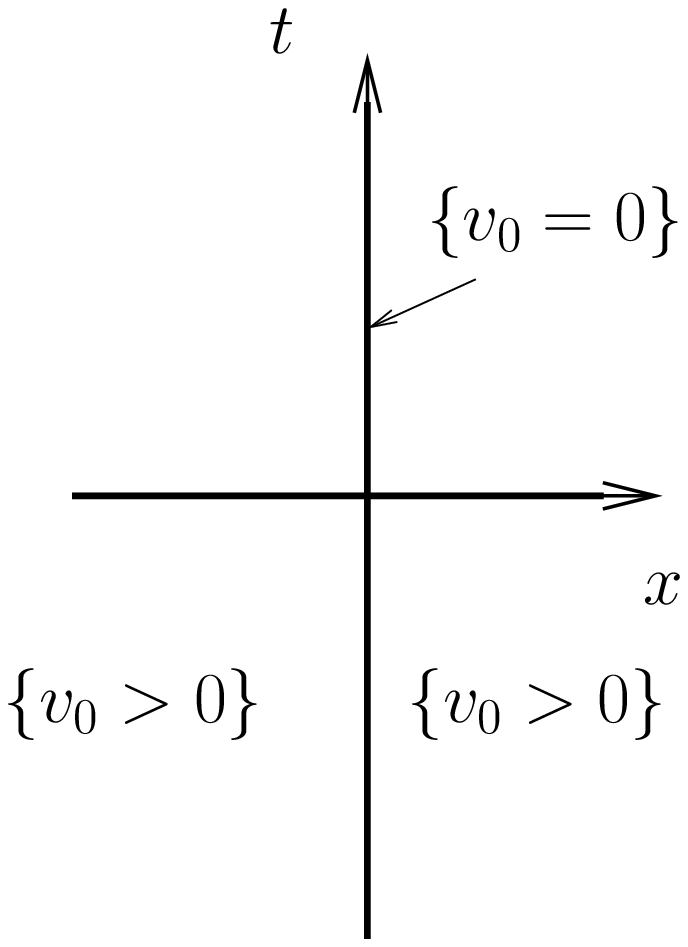}
\end{minipage}
\end{center}
 {\footnotesize \hspace{1.3cm}$v_+$ \hspace{1.9cm} $v_-$ \hspace{1cm} $v_m$, $m \in (-1,0)$\hspace{1.5cm} $v_{-1}$ \hspace{2.1cm} $v_0$}
\caption{\footnotesize Solutions of Theorem \ref{classirr}}
\end{figure}

Similar versions of this theorem are also proved in \cite{CPS}.

Theorem \ref{th:2} is a consequence of Theorem \ref{classirr}. Every blow-up limit satisfies
$u^0_t \le 0$. A more detailed analysis leads to 
\begin{equation}\label{eq:7}
\liminf_{P\to P_0}u_t(P)\le 0\;.
\end{equation}
From the assumption $u_t\ge 0$ we so infer that $u_t = 0$.

We also have an energy criterion to characterise points of the free boundary
\begin{theorem}[Regular and singular points]\label{th:5}
Let $u$ be a solution of (\ref{etoilejean}). Then either ${\mathcal E}(0^-;u) =\sqrt{2}$ or ${\mathcal E}(0^-;u)=\sqrt{2}/2$.
\end{theorem}
 In the first case ({\it i.e.} ${\mathcal E}(0^-;u) =\sqrt{2}$) $P_0$~is
called a {\it singular point}. Otherwise ({\it i.e.} ${\mathcal
  E}(0^-;u)=\sqrt{2}/2$), $P_0$ is a {\it regular point}.\\
\proof
By (\ref{eq:8}) we have ${\mathcal E}(0^-;u)= {\mathcal E}(-1;u^0)$.
Blow-up limits have been classified in Theorem \ref{classirr}. A
simple calculation gives ${\mathcal E}(-1;v_+)= {\mathcal
  E}(-1;v_+)= \sqrt{2}/2$, and ${\mathcal E}(-1;v_m)=\sqrt{2}$ for every
$m\in [-1,0]$.
\finprf
\section{A monotonicity formula for singular points}\label{singulier}
One of the crucial idea of \cite{mono2} can be adapted to the parabolic framework.
\begin{theorem}[Monotonicity formula for singular points]\label{th:6}
Let $u$ be a solution of~(\ref{etoilejean}) and assume that $P_0=0$ is a singular point. For any $m\in [-1,0]$ the function 
\begin{equation}\label{eq:11}
 t \mapsto \Phi^{v_m}(t;u)=\int_{\mathbb R}\frac{1}{t^2}\left(u(x,t)-v_m(x,t)\right)^2 \ G(x,t) \ dx \;, \quad t<0 
\end{equation}
is non-increasing.
\end{theorem}
As a consequence $\lim_{\tau \to 0 ,\; \tau <0} \Phi^{v_m}(t;u):=\Phi^{v_m}(0^-;u)$ is well defined.
\begin{corollary}[Blow-up limit at singular points]\label{cor:7}
Under the assumption of Theorem \ref{th:6}, there exists a real $m\in [-1,0]$
such that any blow-up limit of $u$ at $0$ is equal to $v_{m}$.   
\end{corollary}
\proof
Consider two sequences $(\varepsilon_n)_n$ and $(\tilde\varepsilon_n
)_n$ converging to 0 such that $(u^{\varepsilon_n}_{0})_n$  and
$(u^{\tilde\varepsilon_n}_{0})_n$ respectively converge to two blow-up limits $u^0=v_m$, for some $m \in [0,1]$, and $\tilde u^{0}$. Using the scale-invariance we get
\begin{equation*}
 \Phi^{v_{m}}(-1;u^{\varepsilon_n}_{0})=
 \Phi^{v_{m}}(-\varepsilon_n^2t;u) \quad \mbox{and} \quad \Phi^{v_{m}}(-1;u^{\tilde\varepsilon_n}_{0})=
 \Phi^{v_{m}}(-\tilde\varepsilon_n^2t;u) \;.
\end{equation*}
Passing to the limit in the scale-invariance we obtain
\begin{equation*}
0=\Phi^{v_{m}}(-1;u^{0})=\Phi^{v_{m}}(0^-;u) =\Phi^{v_{m}}(-1;\tilde u^{0})\;.
\end{equation*}
This proves that $\tilde u^{0}=v_{m}=u^{0}$, {\it i.e.} the uniqueness of the blow-up limit.
\finprf

Heuristically the singular points whose blow-up limit is $v_m$ with
$m\in [-1,0)$ have a free boundary with horizontal tangent in the
$(x,t)$-plane: this almost never occurs. On the other hand the singular points with $v_0$ as blow-up limit and the
regular points have a blow-up limit which satisfies $u_t^0=0$. This last
argument can be refined to show that the time derivative of
$u$ is continuous at such points, and consequently for almost every time. This proves Theorem \ref{th:1}.
\medskip\\
\noindent {\bf \large Acknowledgements}\\
We thank D. Lamberton and
B. Lapeyre for their comments and suggestions. This study has been
partially supported by the ACI NIM ``EDP et finance, \# 2003-83''.

\bigskip\noindent{\small \copyright~2004 by the authors. This paper may
be reproduced, in its entirety, for non-commercial purposes.}
\bibliographystyle{siam}
\bibliography{bdm-note-hal}
\end{document}